%% file: agt-1-8.tex
\documentclass{gtart}

\input agtout

\volumenumber{1}
\volumeyear{2001}
\papernumber{8}
\published{21 March 2001}
\pagenumbers{153}{172}
\received{15 December 2000}
\revised{13 February 2001}
\accepted{13 February 2001}

\usepackage{amsmath,amssymb,amsthm}

\newcommand{\SL}{\mbox{\rm SL}}
\newcommand{\SO}{\mbox{\rm SO}}
\newcommand{\Int}{\mbox{\rm Int}}
\newtheorem{thm}{Theorem}
\newtheorem{lem}[thm]{Lemma}
\newtheorem{prop}[thm]{Proposition}
\newtheorem{cor}[thm]{Corollary}
\theoremstyle{definition}
\newtheorem*{defn}{Definition}
\newtheorem*{rem}{Remark}
\newtheorem*{ack}{Acknowledgements}

\begin{document}
\title{Symplectic fillability of tight
contact structures\\on torus bundles}
\authors{Fan Ding\\Hansj\"org Geiges}
\asciiauthors{Fan Ding\\Hansjorg Geiges}
\coverauthors{Fan Ding\\Hansj\noexpand\"org Geiges}
\shorttitle{Fillability of tight contact structures}
\address{Department of Mathematics, Peking University\\Beijing 100871,
P. R. China}
\secondaddress{Mathematisch Instituut, Universiteit Leiden\\Postbus
9512, 2300 RA Leiden, Netherlands}
\email{dingfan@sxx0.math.pku.edu.cn}
\secondemail{geiges@math.leidenuniv.nl}
\asciiaddress{Department of Mathematics, Peking University\\Beijing 100871,
P. R. China\\Mathematisch Instituut, Universiteit Leiden\\Postbus
9512, 2300 RA Leiden, Netherlands}
\asciiemail{dingfan@sxx0.math.pku.edu.cn, geiges@math.leidenuniv.nl}

\asciiabstract{We study weak versus strong symplectic fillability of
some tight contact structures on torus bundles over the circle. In
particular, we prove that almost all of these tight contact structures
are weakly, but not strongly symplectically fillable. For the 3-torus
this theorem was established by Eliashberg.}

\begin{abstract}
We study weak versus strong symplectic
fillability of some tight contact structures on torus bundles
over the circle. In particular, we prove that almost all of these tight contact
structures are weakly, but not strongly symplectically fillable. For the
$3$--torus this theorem was established by Eliashberg.
\end{abstract}

\maketitle

\section{Introduction}

A coorientable $2$--plane field on an oriented $3$--manifold $M$
is called a {\it (positive) contact structure} if, for any $1$--form
$\alpha$ defining $\xi$ as $\xi=\ker\alpha$, the $3$--form $\alpha\wedge
d\alpha$ is a (positive) volume form on~$M$. Notice that the sign of
$\alpha\wedge d\alpha$ only depends on~$\xi$, not on the choice
of~$\alpha$. In this paper, our contact structures are always understood
to be positive. We do not consider non-coorientable contact structures
(where the corresponding $\alpha$ only exists locally).

There are various notions of fillability of contact structures,
see the survey~\cite{etny98}. The two that we are concerned with in the present
paper are weak and strong symplectic fillability.
Given a $4$--dimensional symplectic manifold $(W,\omega )$, we orient it
by regarding $\omega^2$ as a {\it positive} volume form. If $W$ has boundary
$\partial W\! ,$ an orientation of $\partial W$ is defined by the
volume form $i_Y\omega^2$, where $Y$ is any vector field defined along
the boundary and pointing outwards.
Recall that the condition for a vector field $X$ on a symplectic
manifold $(W,\omega )$ to be a Liouville vector field is that
$\mathcal{L}_X\omega =\omega$. By the Cartan formula for the Lie derivative
this may be rewritten as $d(i_X\omega )=\omega$, and this easily implies
that $i_X\omega$ defines a contact structure on any hypersurface transverse
to~$X$.

\begin{defn}
(a)\qua A contact manifold $(M,\xi )$ is {\it weakly symplectically
fillable} if $M$ is the boundary of a symplectic manifold $(W,\omega )$
with $\omega |\xi$ non\-degenerate along $\partial W=M$, and the orientations
on $M$ induced by $W$ and $\xi$ agree.

(b)\qua A contact manifold $(M,\xi )$ is {\it strongly symplectically
fillable} if $M$ is the boundary of a symplectic manifold $(W,\omega)$
admitting a Liouville vector field $X$ near the boundary $\partial W=M$,
pointing outwards along $\partial W\! ,$
and such that $\xi =\ker (i_X\omega | M)$.
\end{defn}

Recall that a contact structure $\xi$ on a $3$--manifold $M$ is called
{\it overtwisted} if there is an embedded 2--disc $D\hookrightarrow M$
such that $\partial D$ is tangent to~$\xi$, but $D$ is transverse to
$\xi$ along~$\partial D$; such a disc is called an {\it overtwisted disc}.
If no such $D$ exists, then $\xi$ is called {\it tight}. Any weakly
symplectically fillable contact structure is tight, as was shown by
Eliashberg and Gromov, cf.~\cite{etny98}.

Clearly strong symplectic fillability implies weak symplectic
fillability. The converse
was shown to be false by Eliashberg~\cite{elia96}. On the $3$--torus
$T^3={\mathbf R}^3/{\mathbf Z}^3$ with coordinates $(x,y,t)$ and orientation
given by $dx\wedge dy\wedge dt$, consider, for non-negative integers~$n$,
the contact structures $\zeta_n$,
defined by
\[ \cos (2\pi (n+1)t)\, dx-\sin (2\pi (n+1)t)\, dy =0.\]
The $\zeta_n$, $n\in {\mathbf N}_0$, are pairwise nondiffeomorphic and
constitute a complete list, up to diffeomorphism, of the tight contact
structures on~$T^3$.

As observed by Giroux~\cite{giro94}, the $\zeta_n$ are all
weakly symplectically fillable. Eliashberg~\cite{elia96} showed that
$\zeta_n$ is strongly symplectically fillable if and only if $n=0$.
Our aim in the present paper is to prove an analogous result for more
general $T^2$--bundles over~$S^1$.

We begin with a description of these torus bundles.
For each matrix $A\in \SL_2({\mathbf Z})$, let $T_A$ denote the quotient of
$T^2\times {\mathbf R}=({\mathbf R}^2/{\mathbf Z}^2)
\times {\mathbf R}$ with coordinates
$({\mathbf x},t)=\bigl( \left( \begin{array}{c} x \\ y \end{array} \right)
,t\bigr)$ 
by the transformation $({\mathbf x},t)\to (A{\mathbf x},t+1)$.
We orient $T_A$ by the
3--form $dx\wedge dy\wedge dt$. The $T^2$--bundle $T_A$ over
$S^1$ depends, up to diffeomorphism, only
on the conjugacy class of $A$ in $\SL_2({\mathbf Z})$.
If $A$ is of the form $\left ( \begin{array} {cc} 1 & 0 \\ k & 1  \end{array}
\right )$, $k\in {\mathbf Z}$, then we denote the corresponding
manifold $T_A$ by $T(k)$.

Let $\varphi\co {\mathbf R}\rightarrow {\mathbf R}$ be a smooth function
whose derivative is strictly positive. The equation
\[ \cos\varphi(t)\, dx- \sin\varphi(t)\, dy=0,\;\;
(x,y,t)\in{\mathbf R}^3,\]
defines a contact structure on ${\mathbf  R}^3$ which we denote by
$\widetilde{\zeta}(\varphi)$. For each $\theta\in {\mathbf R}$ let
$\Delta_{\theta}$ denote the ray
\[ \bigl\{ \left( \begin{array}{c} s\sin\theta
\\ s\cos\theta\end{array} \right)\co s\ge 0\bigr\}\subset {\mathbf R}^2.\]
If $A(\Delta_{\varphi(t)})=\Delta_{\varphi(t+1)}$ for all $t\in
{\mathbf R}$, then the contact structure $\widetilde{\zeta}(\varphi)$ on
${\mathbf R}^3$ is invariant under the action of the deck transformation
group of $T_A$ and thus descends to a contact
structure on $T_A$ which we denote by $\zeta (\varphi)$.

By~\cite{giro99}, for each
non-negative integer $n$ there exists a smooth function
$\varphi\co {\mathbf R}\to {\mathbf R}$ with strictly positive
derivative, satisfying $A(\Delta_{\varphi(t)})
=\Delta_{\varphi(t+1)}$ for all $t\in {\mathbf R}$ and 
\[ 2n\pi<\sup_{t\in {\mathbf R}}\bigl(\varphi(t+1)-\varphi(t)
\bigr)\le 2(n+1)\pi.\]
Up to fibre preserving isotopy, the contact structure $\zeta (\varphi)$ on
$T_A$ depends only on~$n$. Thus we denote this contact structure simply by
$\zeta_n$.  In \cite{giro99} it was shown that the $\zeta_n$ are tight and
pairwise nondiffeomorphic.

The main result of the present paper is the following.

\begin{thm}
For each $A\in\SL_2({\mathbf Z})$ and $n\in {\mathbf N}_0$ (non-negative
integers), the contact
manifold $(T_A,\zeta_n)$ is weakly symplectically fillable. There
exists $n(A)\in {\mathbf N}_0$ such that $(T_A,\zeta_n)$ is not
strongly symplectically fillable for $n>n(A)$.
\end{thm}

Combining this with the classification of tight contact structures
on $T^3_A$ due to Giroux~\cite{giro00} and Honda~\cite{hond00II},
we obtain the following corollary.

\begin{cor}
If $A\in SL_2({\mathbf Z})$ with $\mbox{\rm trace}(A)\neq -2$,
then there are only finitely many strongly symplectically fillable contact 
structures on $T_A$ up to diffeomorphism.
\end{cor}

\begin{proof}
By \cite[Thm.~1.3]{giro00}, cf.~\cite[Thm.~6]{giro99}, $T_A$ admits only
finitely many tight contact structures next to the $\zeta_n$, provided
that $\mbox{\rm trace}(A)\neq -2$. In the case $\mbox{\rm trace}(A)=-2$
there is a further infinite family of tight contact structures.
\end{proof}

Additional results for the $T(k)$ are given in Corollary~\ref{cor:Tkfill}
and Proposition~\ref{prop:Tkfill}.
\section{Preliminaries}
In this section, we review some basic concepts and results needed later.
See \cite{giro91} and \cite{hond00I} for details.

Let $(M,\xi)$ be a contact 3--manifold.
Let $\Sigma$ be an orientable surface embedded in $(M,\xi)$. Let $Y$ be the
vector field on $\Sigma$ defined by the equation $i_Y\Omega=\alpha|_{\Sigma}$,
where $\alpha$ is a global 1--form which defines $\xi$, and
$\Omega$ is an area form on $\Sigma$. The {\it characteristic foliation}
$\xi|_{\Sigma}$ on $\Sigma$ induced by $\xi$ is the singular
foliation represented by $Y$.

A vector field on $M$ is called {\it contact} if its flow preserves $\xi$.
A closed orientable surface $\Sigma$ embedded in $(M,\xi)$ is called
{\it convex} if there is a contact vector field $X$ transverse
to $\Sigma$. This contact vector field $X$ allows us to find a
vertically invariant neighbourhood $\Sigma\times{\mathbf R}\subset M$ of
$\Sigma$, where $\Sigma$ is identified with $\Sigma\times\{ 0\}$.
The {\it dividing set
$\Gamma_{\Sigma}$ for} $X$ is the set of points $x\in\Sigma$ where
$X(x)\in\xi(x)$. This dividing set
$\Gamma_{\Sigma}$ is a disjoint union of simple closed curves which
are transverse to the characteristic foliation $\xi|_{\Sigma}$. The
isotopy type of $\Gamma_{\Sigma}$ is independent of the choice of~$X$. Hence we
will slightly abuse notation and call $\Gamma_{\Sigma}$ {\it the dividing set
of} $\Sigma$. Denote the number of connected components of $\Gamma_{\Sigma}$
by $\#\Gamma_{\Sigma}$.

Let $T$ be a convex torus in a tight contact 3--manifold. Then
the dividing set $\Gamma_T$ consists of an even number $\#\Gamma_T$
of parallel essential curves. Fix an identification of
$T$ with ${\mathbf R}^2/{\mathbf Z}^2$. After a diffeomorphism isotopic
to the identity, we may assume that the dividing curves are linear.
We call the slope of the dividing curves the {\it slope} of the convex
torus $T$ and denote it by~$s(T)$.

Let $V$ be a solid torus.  A specified homeomorphism $h\co S^1\times
D^2\to V$ is called a {\it framing} of $V$.
Fixing such a framing, we identify $\partial V$ with
$T^2={\mathbf R}^2/{\mathbf Z}^2$ by letting
$\{ \left ( \begin{array}{c} t \\ 0 \end{array} \right ) :0\le t\le 1\}$
correspond to the meridian of the solid torus $V$, and 
$\{ \left (  \begin{array}{c} 0 \\ t \end{array} \right ) :0\le t\le 1\}$
correspond to the longitude determined by the framing. With these
identifications, the meridian has slope~$0$, the longitude slope~$\infty$.

The following proposition will prove useful later; see
\cite[Thm.~8.2]{kand97} and \cite[Prop.~4.3]{hond00I}.

\begin{prop}
\label{hon}
For any integer $k$ (including~$0$) there exists a unique (up to
isotopy fixed at the boundary) tight contact
structure on $S^1\times D^2$ with a fixed convex boundary with
$\#\Gamma_{\partial(S^1\times D^2)}=2$ and slope $s(\partial(S^1\times
D^2))=1/k$.
\end{prop}

\begin{rem}
In \cite{hond00I} this is stated for integers~$k$ of a particular sign only.
But if slope $1/k$ can be realised (uniquely), then slope $1/(k+l)$ can be
realised (uniquely) for any integer~$l$, since the two slopes are
related to each other by an $l$--fold Dehn twist along the meridian,
which extends to a diffeomorphism of the solid torus.
\end{rem}

Let $F$ be a closed orientable surface and
$\cal F$ be a singular foliation on $F$. Let $\Gamma$ be a disjoint union of
simple closed curves embedded in $F$ which are transverse to $\cal F$.
Let $F_{\Gamma}$ denote the compact surface with boundary obtained by 
cutting $F$ along $\Gamma$.
We say that $\Gamma$ \it divides \rm $\cal F$ if there 
is a vector field $Y$ on $F_{\Gamma}$ such that

\begin{itemize}
\item $Y$ represents the singular foliation on $F_{\Gamma}$
induced by $\cal F$;
\item $L_Y\Omega>0$ for an area form $\Omega$ on $F_{\Gamma}$;
\item $Y$ goes outward along $\partial F_{\Gamma}$.
\end{itemize}
Here is an important result concerning convex surfaces:

\begin{prop}[Giroux \mbox{\cite[Prop.~II.3.6]{giro91}}]
\label{gir1}
Let $\Sigma$ be a closed convex surface in a
contact $3$--manifold $(M,\xi)$ with contact vector field $X$ and dividing set
$\Gamma$ for~$X$. If $\cal F$ is a singular foliation on $\Sigma$ divided by
$\Gamma$, then there is an isotopy $\phi_s$, $s\in [0,1]$, of $\Sigma$ such
that $\phi_0={\rm id}_{\Sigma},\ \xi|_{\phi_1(\Sigma)}=\phi_1(\cal F)$ and
$\phi_s(\Sigma)$ is transverse to $X$ for each~$s$. 
\end{prop}

\noindent
We state two other results, essentially due to Giroux,
which will be used in Section~4.

\begin{prop}
\label{gir2}
Let $F$ be a closed orientable surface
embedded in a contact $3$--manifold $(M,\xi)$. If $\Gamma$ divides
$\xi|_{F}$, then $F$ is convex with dividing set $\Gamma$.
\end{prop}

\noindent
This proposition is a consequence of~\cite{giro91}, Propositions I.3.4 and
II.1.2(b).

\begin{prop}
\label{gir3}
Let $F$ be a compact orientable surface with 
boundary. Let $\xi_0,\xi_1$ be two contact structures on $F\times
{\mathbf R}$. Let $U$ be a collar neighbourhood of $\partial F$ in $F$.
Assume that $\xi_0|_{F\times\{ 0\}}$ coincides with $\xi_1|_{F\times\{ 0\}}$
and $\xi_0=\xi_1$ on $U\times (-\epsilon,\epsilon)$, where $\epsilon$
is a positive real number. Then there exist a neighbourhood $V$ of $F\times
\{ 0\}$ in $F\times {\mathbf R}$ and a contact embedding $f\co (V,\xi_0)\to
(F\times {\mathbf R},\xi_1)$
such that $f$ is the identity on $V\cap ( U\times (-\epsilon,\epsilon ))$
and $f(F\times\{ 0\})=F\times \{ 0\}$.
\end{prop}

\noindent
The proof of this proposition is similar to that of \cite{giro91},
Proposition II.1.2(b).
\section{Contact surgery}
A smooth knot $K\co S^1\to M$ in a contact $3$--manifold $(M,\xi)$ is called
{\it Legendrian} if its tangent vectors all lie in $\xi$.
Any diffeomorphism between Legendrian knots extends to a contactomorphism
(i.e.\ a diffeomorphism preserving contact structures) on
some neighbourhoods of the knots. A Legendrian knot $K$ comes 
equipped with a canonical framing of its normal bundle, which is
induced by any vector field transverse to $\xi$, or equivalently, by
a vector field in $\xi|K$ transverse to $K$. We call this the
{\it contact framing} of~$K$. There is
a canonical bijection from (normal) framings of $K$ to the integers
${\mathbf Z}=\pi_1(\SO (2))$, given by identifying the contact framing
with $0\in\mathbf{Z}$ and counting right-handed twists positively.
Note that for any nullhomologous
Legendrian knot $K$, the linking number of $K$ with its push-off determined
by framing $k$ is $tb(K)+k$, where $tb(K)$ is the {\it Thurston-Bennequin
invariant} of~$K$.

{\it Rational surgery} on $K$ with coefficient $r=p/q\in {\mathbf Q}\cup
\{ \infty\}$ (with $p,q$ coprime) is defined as follows: Denote a tubular
neighbourhood of K
(diffeomorphic to a solid torus) by $\nu K$. Let $(\mu,\lambda)$ be a
positively oriented basis for $H_1(\partial\nu K;{\mathbf Z})\cong
{\mathbf Z}\oplus {\mathbf Z}$, where $\lambda$ is determined up to sign as
the class of a parallel copy of $K$ determined by the contact framing, and
$\mu$ is determined  by a suitably oriented meridian (i.e.\ a
nullhomologous circle in $\nu K$), cf.~\cite[p.~672]{gomp98}.
We obtain a new manifold $M'$ by cutting
$\nu K$ out of $M$ and regluing it by a diffeomorphism of $\partial(\nu K)$
sending $\mu$ to $p\mu+q\lambda$. This procedure determines $M'$ up to
orientation-preserving diffeomorphism. 

Consider $N={\mathbf R}^2\times ({\mathbf R}/{\mathbf Z})$ with coordinates
$(x,y,z)$ and contact structure $\zeta$ defined by
\[ \cos (2\pi z)\, dx-\sin (2\pi z)\, dy =0.\]
For each $\delta >0$, let
\[ N_{\delta}= \{ (x,y,z)\in N\co x^2+y^2\le {\delta}^2\}.\]
We identify $\partial N_{\delta}$ with ${\mathbf R}^2/{\mathbf Z}^2$,
using the contact framing, and
write $(\mu ,\lambda )$ for a positively oriented basis for
$H_1(\partial N_{\delta};{\mathbf Z})\cong {\mathbf Z}\oplus {\mathbf Z}$,
with $\mu$ corresponding to a meridian and $\lambda$ to a longitude
determined by this framing. A possible representative of $\lambda$ would be
\[ \{ (\delta\sin (2\pi z),\delta\cos (2\pi z),z)\co z\in
{\mathbf R}/{\mathbf Z}\}.\]

Note that the vector field $x\frac{\partial}{\partial x}+y
\frac{\partial}{\partial y}$ is a contact vector field for~$\zeta$ which
is transverse to $\partial N_{\delta}$, with dividing set
\[ \Gamma_{\partial N_{\delta}}=\{ (\pm\delta\sin (2\pi z),\pm\delta\cos
(2\pi z),z)\co z\in {\mathbf R}/{\mathbf Z}\} .\]
Thus for each ${\delta}>0$ the torus $\partial N_{\delta}$ is a convex surface
with $\#\Gamma_{\partial N_{\delta}}=2$ and $s(\partial N_{\delta})=\infty$.

Let $K$ be a Legendrian knot in a contact manifold $(M,\xi)$. Write
\[ C=\{ (x,y,z)\in N\co x=y=0\} \]
for the spine of~$N$. Then there is a contact embedding $f\co (N_2,\zeta)\to
(M,\xi)$ such that $f(C)=K$. We want to construct a contact structure
$\xi '$ on the manifold $M'$ obtained from $M$ by rational surgery on
$K$ with coefficient~$r=p/q\in {\mathbf Q}\cup \{\infty\}$, where
we only consider $r\neq 0$. Let
\[ P=\{(x,y,z)\in N\co 1\le x^2+y^2\le 4\}=N_2\setminus\Int\, N_1.\]
Let $g\co P\to P$ be an orientation-preserving diffeomorphism
sending $\partial N_{\delta}$ to $\partial N_{\delta}$, ${\delta}=1,2$, and
$\mu$ to $p\mu +q\lambda$. The fact that $p\neq 0$ implies that
$(g_*)^{-1}(\zeta)$ is a contact structure on $P$ with respect to
which $\partial N_{\delta}$ is a convex torus of
non-zero slope. By~\cite[Thm.~2.3]{hond00I}, which gives
an enumeration of tight contact structures on the solid torus with
convex boundary as in our situation (and in particular shows this
set of contact structures to be non-empty), the contact structure
$(g_*)^{-1}(\zeta )$ on $P$ can be extended to a tight
contact structure $\zeta '$ on~$N_2$. Define
\[ M' =(M-f(N_1))\cup N_2/\sim ,\]
where $x\in P\subset N_2$ is identified with $f(g(x))\in M$.
Topologically, $M'$ is obtained from $M$ by rational surgery on $K$
with coefficient~$r$. It inherits a contact structure~$\xi '$ from
$(M,\xi )$ and $(N_2,\zeta ')$. We say that $(M',\xi ')$ is obtained from
$(M,\xi )$ by {\it contact $r$--surgery on}~$K$.

\begin{rem}
(1) In this construction the assumption $r\neq 0$ is essential. The $g\co
P\rightarrow P$ corresponding to $p=0$, $q=1$ leads to a contact structure
$(g_*)^{-1}(\zeta )$ on $P$ whose extension to $N_2$ (if such
exists) is overtwisted;
the overtwisted disc being given essentially by a meridianal disc
in the solid torus~$N_1$.

(2) It is not clear {\it a priori} that $(M', \xi ')$ is tight, even
if $(M,\xi )$ was. In the
application of this construction (Proposition~\ref{prop:AA0}, in particular)
we deal with a situation where one knows two tight
contact manifolds $(M,\xi )$ and $(M',\xi ')$
to be contactomorphic outside certain solid tori, and we can conclude
there that one is obtained from the other by contact surgery as described.

(3) By analysing the framing conditions in the surgery theorems of
\cite{elia90} and \cite{wein91}, cf.~\cite[Thm.~1.3]{gomp98} and~\cite{etho},
one sees that contact $(-1)$--surgery corresponds to a symplectic
handlebody construction. In particular, if
$(M',\xi')$ is obtained from a closed contact manifold $(M,\xi)$
by contact $(-1)$--surgery and $(M,\xi)$ is strongly symplectically fillable,
then $(M',\xi')$ is also strongly symplectically fillable. Given a
Legendrian knot, one can add left-twists to its contact framing
by performing a suitable isotopy (non-contact and $C^0$--small).
That way one can realise
topological surgeries with negative integer framing (relative to a given
contact framing) as `handlebody' surgeries. Adding positive twists is not,
in general, possible, unless the contact structure is overtwisted. We
are mostly concerned with contact $(1/k)$--surgeries, $k\in
{\mathbf Z}\setminus\{ 0\}$, which do not correspond to a handlebody
construction unless $1/k=-1$.
\end{rem}

\begin{prop}
\label{ind}
If $r=1/k$, where $k$ is an integer, then, up to contactomorphism, the
contact manifold $(M',\xi')$ depends only on $r$ (and $(M,\xi )$ and
$K$, of course). That is, it is independent of
the choices of $f,g$ and~$\zeta '$.
\end{prop}

\begin{proof}
The scaling map $(x,y,z)\mapsto (sx,sy,z)$ defines
a contactomorphism $(N_{\delta},\zeta )\rightarrow (N_{s\delta},
\zeta )$. Hence, given two contact embeddings $f_i\co (N_2,\zeta )
\rightarrow (M,\xi )$, $i=1,2$, we can compare either with a third such
embedding that maps $N_2$ into the interior of $f_i(N_1)$.

We may therefore assume that $K=C\subset N$, the contact embedding
$f_1$ is the inclusion map $N_2\subset N$, and the contact embedding
$f=f_2$ sends $N_2$ into the interior $\mbox{\rm Int}\, N_1$ of~$N_1$.

Note that if $r=1/k$, then the diffeomorphism $g$ may be assumed to have the
following effect on $\mu$ and $\lambda$, since $g|\partial N_{\delta}$
is determined up to isotopy by its action on homology, corresponding
to an element of $\SL_2({\mathbf Z})$:
\[ \mu\longmapsto \mu +k\lambda ,\;\; \lambda\longmapsto
\lambda -l(\mu +k\lambda ) ,\]
where $l$ is some integer. (Different choices of $l$ correspond to Dehn
twists along a meridian of the solid torus that is glued back; these Dehn
twists extend to diffeomorphisms of the solid torus and hence have no
topological effect.) Then $g^{-1}$ sends
$\lambda$ to $l\mu +\lambda$.
This implies that as we pull back $\zeta$ to $(g_*)^{-1}\zeta$, we obtain a
contact structure on $N_2\setminus \Int (N_1)$ with $\#\Gamma_{\partial N_2}
=2$ and $s(\partial N_2)=1/l$. So by Proposition~\ref{hon} the extension of
$(g_*)^{-1}\zeta$ to a tight contact structure $\zeta '$ on the copy
of $N_2$ to be glued back is unique.

Let $(M_1',\xi_1')$ be the contact manifold obtained from $N_2$ by
contact $r$--surgery along $C\subset N_2$ using the inclusion $N_2
\subset N_2$, and let $(M_2',\xi_2')$ be the contact manifold obtained
similarly using the contact embedding $f\co N_2\rightarrow \mbox{\rm Int}
\, N_1\subset N_2$.

By what we have just observed, the tight contact structure $\xi_1'$ is uniquely
determined by the fact that it coincides with $\zeta$ near $\partial N_2
=\partial M_1'$. We also know that $\xi_2'$ coincides with $\zeta$
outside $f(N_1)$, and the manifolds $M_1'$ and $M_2'$ are diffeomorphic
under a diffeomorphism that is the identity near $\partial M_1'=\partial
M_2'$.

By the definition of contact $r$--surgery, the contact manifold
\[ (M_2'\setminus (N_2\setminus f(N_2)), \xi_2')\]
is tight. It now suffices to show that $(M_2',\xi_2')$ is tight, because
we then know that it is completely determined by its boundary data,
which coincide with those of $(M_1',\xi_1')$.

Recall that if a contact structure $\xi$ on a manifold $M$ is written as the
kernel of a $1$--form $\alpha$, there is a one-to-one correspondence
between contact vector fields $X$ and functions on $M$ given by
$X\mapsto \alpha (X)$, cf.~\cite{lima87}. The function $H=\alpha (X)$
is called the {\it Hamiltonian function} corresponding to~$X$.

So the contact vector field $X=-(x\frac{\partial}{\partial x}+
y\frac{\partial}{\partial y})$ for $\zeta$ on $N_2$ corresponds in this way
to some Hamiltonian function. By multiplying this function with a bump
function that is identically $1$ on $N_1$ and identically zero near $\partial
N_2$ we can construct a contactomorphism $N_2\rightarrow N_2$ that is
the identity near $\partial N_2$ and sends $N_1$ into $N_{\delta}$ for any
given $\delta >0$.  By precomposing $f$ with such a diffeomorphism, we may
assume that
\[ f(N_1)\subset\mbox{\rm Int}\,  N_{\delta},\;\;
N_{\delta}\subset\mbox{\rm Int}\, f(N_2)\]
for a suitable $\delta >0$.

By multiplying the Hamiltonian function of $X$ with a bump function that
is identically $0$ on $f(N_1)$ and identically $1$ outside $N_{\delta}$,
we get a Hamiltonian function defined also on~$M_2'$
whose contact flow will ultimately move $N_2$ into $f(N_2)$. So this will
define a contact embedding
\[ (M_2',\xi_2')\hookrightarrow (M_2'\setminus (N_2\setminus f(N_2)),
\xi_2').\]
This completes the proof of the proposition.
\end{proof}

\begin{prop}
\label{inv}
If $(M',\xi')$ is obtained from $(M,\xi)$ by contact $(1/k)$--surgery,
then $(M,\xi)$ is obtained from $(M',\xi')$ by contact $(-1/k)$--surgery.
\end{prop}

\begin{proof}
By the preceding proposition it suffices to consider the following
situ\-ation: Let $(M',\xi ')$ be the manifold obtained from $(N,\zeta )$
by contact $(1/k)$--surgery along $C\subset N$, using the inclusion $N_2
\subset N$. Let $(M,\xi )$ be the manifold obtained from $(M',\xi ')$
by contact $(-1/k)$--surgery along a spine of the
solid torus $N_1$ that was attached to $N\setminus N_1$ to form~$M'$.
We want to show that $(M,\xi )$ is contactomorphic to $(N,\zeta )$.

We can obtain $(M',\xi ')$ by gluing $N_2$ to $N\setminus N_1$ using the
attaching map $g\co P\rightarrow P$ described by
\[ \mu\longmapsto \mu +k\lambda ,\;\; \lambda\longmapsto
\lambda ,\]
and then extending $(g_*)^{-1}\zeta$ over $N_2$ to a unique tight contact
structure~$\zeta '$. We observed in the proof of the preceding proposition
that the torus $\partial N_1$ in $(N_2,\zeta ')$ is a convex
surface with $\# \Gamma_{\partial N_1}=2$ and $s(\partial N_1)=\infty$.

By Propositions \ref{hon} and \ref{gir1} and arguments similar to those
in the preceding proof, we can find a contact embedding
$(N_2,\zeta )\hookrightarrow (N_2,\zeta ')$ isotopic to the identity
and sending $P$ into~$P$. Now perform the
$(-1/k)$--contact
surgery on $(M',\xi ')$ using this embedding (composed with~$g$),
and call the resulting contact manifold $(M,\xi )$. The gluing for this
surgery may be described by
\[ \mu\longmapsto \mu -k\lambda ,\;\; \lambda\longmapsto
\lambda .\]
It is a straightforward check that the topological effect of this second
surgery is to cancel the first surgery,
because the composition of these maps sends $\mu$ to $\mu$ (in fact,
it is the identity map).
A further application of Proposition~\ref{hon} shows that $(M,\xi )$
is indeed contactomorphic to~$(N,\zeta )$.
\end{proof}

Let $(M',\xi ')$ be obtained from $(M,\xi )$ by contact $(1/n)$--surgery
on a Legendrian knot~$K$, where $n>1$. Let $(M'',\xi '')$ be obtained from
$(M,\xi )$ by contact $(1/(n-1))$--surgery on the same knot~$K$. By
the same methods as in the proof of the preceding proposition one sees
that $(M',\xi ')$ can be obtained from $(M'',\xi '')$ by contact
$(+1)$--surgery. Similarly, contact $(1/n)$--surgery with $n<-1$ can be
realised as a contact $(1/(n+1))$--surgery followed by a contact
$(-1)$--surgery. Thus, by induction we have:

\begin{prop}
If $(M',\xi ')$ is obtained from $(M,\xi )$ by contact $(1/n)$--surgery,
$n\neq 0$, it may also be obtained by $|n|$ times contact
$\varepsilon$--surgery, where $\varepsilon =n/|n|=\mbox{\rm sign}(n)
\in\{-1,1\}$.
\end{prop}

Combined with remark (3) above, the two preceding propositions yield the
following result.

\begin{prop}
\label{prop:fillsurgery}
Let $(M,\xi)$ be a closed contact $3$--manifold.
Let $n$ be a positive integer.
If $(M',\xi')$ is obtained from $(M,\xi)$ by contact $(-1/n)$--surgery
and $(M,\xi)$ is strongly symplectically fillable,
then $(M',\xi')$ is strongly symplectically fillable.
If $(M',\xi')$ is obtained from $(M,\xi)$ by contact $(1/n)$--surgery
and $(M,\xi)$ is not strongly symplectically
fillable, then $(M',\xi')$ is not strongly symplectically
fillable.
\end{prop}
\section{Proof of the main result}
The key step in the proof of Theorem~1 is contained in the following
proposition.

\begin{prop}
\label{prop:AA0}
Let $A_0\in\SL_2({\mathbf Z})$, let $E_k=\left( \begin{array}{cc}
1&0\\k&1\end{array}\right)$, $k\in {\mathbf Z}\setminus\{ 0\}$, and let $n_0$
be a positive integer. Then contact $(-1/k)$--surgery on $(T_{A_0},
\zeta_{n_0})$ yields $(T_A,\zeta_n)$, where $A=E_kA_0$ and
$n\in\{ n_0,n_0-1\}$
for $k>0$, or $n\in\{ n_0,n_0+1\}$ for $k<0$. If $A_0$ is of type $E_l$,
then $n$ is determined explicitly as follows.

\begin{center}
\begin{tabular}{c|l|c}
$k$ & $A_0$ & $n$\\[.5mm] \hline
$>0$ & $E_l$, $l<-k$ or $l\geq 0$ & $n_0$\\[.5mm]
     & $E_l$, $-k\leq l<0$ & $n_0-1$\\[.5mm]
$<0$ & $E_l$, $l\geq -k$ or $l<0$ & $n_0$\\[.5mm]
     & $E_l$, $0\leq l<-k$ & $n_0+1$\\
\end{tabular}
\end{center}
\end{prop}

\begin{rem}
An exact determination of the value of $n$ corresponding to any given
$A_0$ is feasible and would allow an estimate on the bound $n(A)$ in
Theorem~1.
\end{rem}

Using this proposition, we can formulate a strengthening of Theorem~1
in the case $A=E_k$, $k<0$.

\begin{cor}
\label{cor:Tkfill}
The contact manifold $(T(k),\zeta_n)$ is not strongly symplectically
fillable for $k\leq 0$ and $n\geq 2$.
\end{cor}

\begin{proof}
For $k=0$, that is, $T(0)=T^3$, this is the result of Eliashberg
mentioned in the introduction, which holds true even for $n=1$. By the
preceding proposition, $(T(k),\zeta_{n+1})$, $k<0$, is obtained from
$(T^3,\zeta_n)$ by contact $(-1/k)$--surgery. The result now follows
from Proposition~\ref{prop:fillsurgery}.
\end{proof}

Here is a complementary result.

\begin{prop}
\label{prop:Tkfill}
The contact manifold $(T(k),\zeta_0)$ is strongly symplectically fillable
for all $k\in{\mathbf Z}$.
\end{prop}

\begin{proof}
For $k=0$ this is well-known, see \cite{elia96}. For positive $k$ it is a
consequence of Propositions \ref{prop:fillsurgery} and~\ref{prop:AA0} (which
holds true also for $n_0=0$ and $A_0=E_0$).

For negative $k$ we use a construction analogous to~\cite[Lemma~2.6]{mcdu91}.
Let $\varphi\co {\mathbf R}\rightarrow {\mathbf R}$ be a smooth function
with strictly positive derivative, $E_k(\Delta_{\varphi (t)})=
\Delta_{\varphi (t+1)}$ for all $t\in {\mathbf R}$,
and $\varphi (0)=\pi /2$. Notice that $k<0$ then implies
$0<\varphi (t)<\pi$ for all $t\in {\mathbf R}$. So we may define
$\zeta_0$ as kernel of the  contact form $\beta =dy-\cot \varphi (t)\, dx$
(defined on~$T(k)$).

Projection onto the $x$-- and $t$--coordinate gives $T(k)$ the structure
of a principal $S^1$--bundle $\pi\co T(k)\rightarrow T^2$. Let $L$ be the
associated complex line bundle $T(k)\times_{S^1}{\mathbf C}$, and write
$L_0$ for its zero section. Write $\theta$ for the
angular coordinate and $r$ for the radial coordinate in the
${\mathbf C}$--fibre, so that $\frac{\partial}{\partial y}=
\frac{\partial}{\partial\theta}$ on $T(k)\subset L$. The vector fields
$\frac{\partial}{\partial\theta}$ and $\frac{\partial}{\partial r}$
are defined on $L\setminus L_0$, and $\beta$ extends to an
$S^1$--invariant $1$--form on $L\setminus L_0$ satisfying
$\beta (\frac{\partial}{\partial\theta})=1$ and
$\beta (\frac{\partial}{\partial r})=0$. We then have
\[ d\beta =\pi^*\bigl( \varphi '(t)\csc^2\varphi (t)\, dt\wedge dx\bigr) .\]
Set
\[ \omega = d((r^2+1)\beta )=(r^2+1)d\beta +2r\, dr\wedge d\beta .\]
It is a straightforward check that $\omega$ is a symplectic form defined
on all of~$L$, and that $X=\frac{r^2+1}{2r}\frac{\partial}{\partial r}$
is a Liouville vector field for $\omega$ defined on $L\setminus L_0$,
and $i_X\omega =(r^2+1)\beta$.
So the unit disc bundle $T(k)\times_{S^1}D^2$ gives a strong symplectic
filling of $(T(k),\zeta_0)$ for $k<0$.
\end{proof}

Assuming Proposition~\ref{prop:AA0}, we can now prove the part of Theorem~1
concerned with strong symplectic fillability.

\begin{prop}
For each $A\in\SL_2({\mathbf Z})$ there exists an $n(A)\in {\mathbf N}_0$ such
that $(T_A,\zeta_n)$ is not strongly symplectically fillable for $n>n(A)$.
\end{prop}

\begin{proof}
It is well-known (and easy to prove) that $\SL_2({\mathbf Z})$ is
generated by $E_{-1}$ and $E_1'=\left(\begin{array}{cc}1&1\\0&1\end{array}
\right)$. Moreover, the product $E_{-1}E_1'$ is of order $6$ in
$\SL_2({\mathbf Z})$, which implies that $E_{-1}^{-1}$ and $(E_1')^{-1}$
can be expressed as a pro\-duct in $E_{-1}$ and~$E_1'$. Thus, with $A\in\SL_2(
{\mathbf Z})$ given, we may write it in the form $A=A_1\cdots A_m$ with
$m\in {\mathbf N}$ and $A_i\in\{ E_{-1},E_1'\}$. Set $A_m'=E_0$ and
$A_i'=A_{i+1}\cdots A_m$ for $i=1,\ldots,m-1$, so that $A_{i-1}'=A_iA_i'$.

If $A_i=E_{-1}$, then by Proposition~\ref{prop:AA0} we know that
$(T_{A_{i-1}'},\zeta_n)$ is obtained from $(T_{A_i'},\zeta_{n_0})$ by
contact $(+1)$--surgery, where $n\in\{ n_0,n_0+1\}$ is chosen suitably.

If $A_i=E_1'$, we observe that with $B=\left(\begin{array}{cc}
0&1\\-1&0\end{array}\right)$ we can write
\[ A_{i-1}'=A_iA_i'=B(E_{-1}(B^{-1}A_i'B))B^{-1}.\]
Since conjugate matrices $B_0$ and $B_1=BB_0B^{-1}$ give rise to
contactomorphic torus bundles $(T_{B_0},\zeta_n)$ and $(T_{B_1},\zeta_n)$,
we conclude once again that $(T_{A_{i-1}'},\zeta_n)$ is obtained from
$(T_{A_i'},\zeta_{n_0})$ by contact $(+1)$--surgery for a suitable
$n\in\{ n_0,n_0+1\}$.

By induction, there exists $n(A)\in {\mathbf N}_0$ such that
$(T_A,\zeta_{n+n(A)})$ is obtained from $(T^3,\zeta_n)$ by $m$ times contact
$(+1)$--surgery. Thus, Eliashberg's theorem and
Proposition~\ref{prop:fillsurgery} imply that
$(T_A,\zeta_{n+n(A)})$ is not strongly
symplectically fillable for~$n\geq 1$.
\end{proof}

To prove the part of Theorem~1 concerned with weak symplectic fillability,
we first make the following observation.

\begin{prop}
For each $A\in\SL_2({\mathbf Z})$ there exists a compact symplectic
manifold $(W,\omega )$ such that $T_A$ is the oriented boundary of $W$
and $\omega$ is nondegenerate on each torus fibre of~$T_A$.
\end{prop}

\begin{proof}
Given $A\in\SL_2({\mathbf Z})$, we write it in the form $A=A_1\cdots A_m$ with
each $A_i$ equal to $E_{-1}$ or~$E_1'$. Let $\pi\co S\rightarrow {\mathbf C}
P^1$ be a nodal elliptic surface with a section, without multiple fibres,
and with Euler number (or number of singular fibres) equal to $12d\geq 2m$.
For the existence of such a surface see~\cite[p.~64]{frmo94}. This
surface is algebraic~\cite[p.~34]{frmo94} and thus K\"ahler; in particular we
find a symplectic from $\omega$ on $S$ that restricts to an area form on
each nonsingular fibre (since these are complex submanifolds).

By the arguments in Section~2.3 of \cite{frmo94} we find a simple closed loop
$\gamma$ in ${\mathbf C}P^1$ along which the monodromy of the fibration~$\pi$
equals~$A$. Let $D\subset {\mathbf C}P^1$ be the disc whose oriented
boundary is~$\gamma$. Then $(W=\pi^{-1}(D),\omega )$ is the desired
symplectic manifold.

Here is an alternative and slightly more direct argument:
Observe that $E_{-1}$ and $E_1'$ correspond to positive Dehn twists of~$T^2$.
This implies that there is an orientable Lefschetz fibration $W\rightarrow
D^2$ with generic fibre a torus,
$m$ singular fibres, and monodromy along $\partial D^2$ equal to~$A$,
cf.~\cite[Section~8.2]{gost99}. Such a Lefschetz fibration admits a symplectic
form $\omega$ with the described properties, see~\cite[Thm.~10.2.18]{gost99}.
Since the base of the fibration is~$D^2$, the second homo\-logy group
of the total space is generated by the fundamental class of the fibre
(this remains true in the presence of singular fibres). So the
homological condition in the cited theorem, necessary to apply
Thurston's symplectic fibration construction, is trivially
satisfied.
\end{proof}

Here is the part of Theorem~1 concerned with weak symplectic fillability:

\begin{prop}
For each $A\in\SL_2({\mathbf Z})$ and $n\in {\mathbf N}_0$, the contact
manifold $(T_A,\zeta_n)$ is weakly symplectically fillable.
\end{prop}

\begin{proof}
Represent $\zeta_n$ by $\zeta (\varphi )$, i.e.
\[ \cos\varphi (t)\, dx-\sin\varphi (t)\, dy=0,\]
where $\varphi$ is as described in the introduction. The properties
of $\varphi$ imply that we can find a smooth function $\lambda\co
{\mathbf R}\rightarrow {\mathbf R}^+$ such that the contact $1$--form
\[ \alpha =\lambda (t)(\cos\varphi (t)\, dx-\sin\varphi (t)\, dy)\]
is invariant under the transformation $({\mathbf x},t)\mapsto
(A{\mathbf x},t+1)$ and thus descends to a contact form (which we continue
to denote~$\alpha$) on $T_A$ representing~$\zeta_n$.

Observe that the $1$--form $\alpha_{\varepsilon}=(1-\varepsilon)\,
dt+\varepsilon\alpha$
is a contact form for any $\varepsilon\in (0,1]$, and in view of the well-known
Gray stability theorem~\cite{gray59} it defines a contact structure
equivalent to~$\zeta_n$. For $\varepsilon\searrow 0$ the contact
planes $\ker\alpha_{\varepsilon}$ approach the tangent spaces along the
fibres of~$T_A$. Hence, the symplectic form $\omega$ on $W$ constructed
in the preceding proposition will have the property that
$\omega|\ker\alpha_{\varepsilon}$ is nondegenerate
for $\varepsilon >0$ sufficiently small.
\end{proof}

To complete the proof of our main theorem it remains to prove
Proposition~\ref{prop:AA0}. We only do this for the case $k=1$;
the other cases are analogous.

Let $A_0\in\SL_2({\mathbf Z})$ be given and set $A=E_1A_0$. Let
$\varphi_0\co {\mathbf R}\rightarrow {\mathbf R}$ be a smooth function
with strictly positive derivative, and satisfying $A(\Delta_{\varphi_0(t)})
=\Delta_{\varphi_0(t+1)}$, where $\Delta_{\theta}$ was defined in the
introduction. The non-negative integer $n_0$ determined by
\[ 2n_0\pi < \sup_{t\in {\mathbf R}}\bigl( \varphi_0(t+1)-\varphi_0(t)
\bigr) \leq 2(n_0+1)\pi \]
will be referred to as the {\it twisting} of~$\varphi_0$. Assume in addition
that $\varphi_0(0)=0$.

\begin{lem}
\label{lem:phi}
There is a smooth function $\varphi\co {\mathbf R}\rightarrow {\mathbf R}$
with strictly positive derivative, satisfying $A(\Delta_{\varphi (t)})
=\Delta_{\varphi (t+1)}$, as well as $\varphi (0)=\varphi_0(0)=0$ and
$\varphi (-1)=\varphi_0(-1)$. The twisting $n$ of this function $\varphi$
depends on $A_0$ and $n_0$ as described in Proposition~$\ref{prop:AA0}$.
\end{lem}

\begin{proof}
It is possible to choose the values of $\varphi (t)$ equal to those of
$\varphi_0(t)$ at $t=0$ and $t=-1$ and still satisfy the appropriate
equivariance condition because $\left(\begin{array}{c}0\\1\end{array}
\right) =\left(\begin{array}{c}\sin\varphi_0(0)\\ \cos\varphi_0(0)
\end{array}\right)$ is an eigenvector of $E_1$ with eigenvalue~$1$.

(i)\qua First consider the case that $\left(\begin{array}{c}0\\1\end{array}
\right)$ is an eigenvector of $A_0$ with positive eigenvalue. This is
equivalent to saying that $A_0$ is of type~$E_l$. A straightforward
analysis shows that in this case
\[ \varphi_0(-1)=\left\{ \begin{array}{ll}
-2n_0\pi & \mbox{\rm if}\;\; l<0,\\
-2(n_0+1)\pi & \mbox{\rm if}\;\; l\geq 0.
\end{array}\right. \]
The same analysis applies to $A=E_1A_0=E_{l+1}$. That is, the function
$\varphi$ with the described properties has twisting $n$ determined by
\[ \varphi (-1)=\left\{ \begin{array}{ll}
-2n\pi & \mbox{\rm if}\;\; l+1<0,\\
-2(n+1)\pi & \mbox{\rm if}\;\; l+1\geq 0.
\end{array}\right. \]
Since $\varphi (-1)=\varphi_0(-1)$ by assumption, we have $n=n_0$ for
$l\neq -1$, and $n=n_0-1$ for $l=-1$.

(ii)\qua Now assume that $A_0$ is not of type~$E_l$. Then $\varphi_0(-1)\not\in
2\pi {\mathbf Z}$, and one verifies that the twisting $n_0$ of $\varphi_0$
is determined by
\[ 2n_0\pi <\sup_{t\in [-1,0]}\bigl( \varphi_0(t+1)-\varphi_0(t)\bigr) \leq
2(n_0+1)\pi,\]
cf.~\cite[p.~791]{giro99}. Let $\overline{h}\co S^1\rightarrow S^1$
(with $S^1={\mathbf R}/2\pi {\mathbf Z}$) be the smooth function defined
by $E_1(\Delta_{\theta})=\Delta_{\overline{h}(\theta )}$, and let
$h\co {\mathbf R}\rightarrow {\mathbf R}$ be the lift of $\overline{h}$
with $h(0)=0$. One checks that $h$ is strictly increasing and
$t-\pi /2\leq h(t)\leq t$ for all $t\in
{\mathbf R}$, with equality $h(t)=t$ for $t\in\pi {\mathbf Z}$.

The required function $\varphi$ can be defined by
smoothing the function
\[ \varphi (t)=\left\{ \begin{array}{ll}
\varphi_0(t) & -3/4\leq t\leq 0,\\
h(\varphi_0(t)) & 0\leq t\leq 3/4,
\end{array}\right. \]
at $t=0$, and then extending it to all $t\in {\mathbf R}$
by imposing the appropriate equivariance property.

Since the smoothing is done at $t=0$, and $\left( \begin{array}{c}
\sin\varphi_0(0)\\ \cos\varphi_0(0)\end{array}\right) =
\left( \begin{array}{c}0\\1\end{array}\right)$ is not an eigenvector
of $A_0$ with positive eigenvalue, we can ensure that this does not
lead to a twisting $n$ larger than~$n_0$. The properties of $h$ imply
that this twisting~$n$, determined by
\[ 2n\pi <\sup_{t\in [-1,0]}\bigl( \varphi (t+1)-\varphi (t)\bigr) \leq
2(n+1)\pi,\]
is equal to $n_0$ or $n_0-1$.
\end{proof}

The strategy in the proof of Proposition~\ref{prop:AA0} is now as
follows. Remove a tubular neighbourhood $T^2\times I$ of a torus
fibre in both $(T_{A_0},\zeta_{n_0})$ and $(T_A,\zeta_n)$, and show that
the complements are contactomorphic provided $\zeta_{n_0}$ corresponds
to $\varphi_0$ and $\zeta_n$ to the $\varphi$ constructed in the
preceding lemma. Extend this contactomorphism over a solid torus
inside $T^2\times I$, with complement another solid torus. Finally
show that the unique extensions (as tight contact structures) of the contact
structures $\zeta_{n_0}$ resp.\ $\zeta_n$ over this last solid torus
correspond to a contact $(-1)$--surgery.

The next lemma will be essential for this final extension. Consider
$B={\mathbf R}\times ({\mathbf R}/{\mathbf Z})\times {\mathbf R}$ with
coordinates $(x,y,t)$ and contact structure $\zeta '$ given by
\[ \cos (2\pi t)\, dx-\sin (2\pi t)\, dy=0.\]
For $0< \varepsilon ,\delta <1/4$ let
\[ V=\{ (x,y,t)\in B\co \delta\leq x\leq 1-\delta ,\; -\varepsilon\leq t\leq
\varepsilon\} .\]
This will later be thought of as a tubular neighbourhood in $(T_A,\zeta_n)$
of a Le\-gendrian circle ($t=0$, $x=\mbox{\rm const.}$), which lies
completely inside a torus fibre of~$T_A$.
Identify $\partial V$ (with corners smoothed) with
${\mathbf R}^2/{\mathbf Z}^2$ by using the standard framing of~$V$. This
means that the circles $y=\mbox{\rm const.}$ (oriented positively in
the $(t,x)$--plane) correspond to the first coordinate direction in
${\mathbf R}^2/{\mathbf Z}^2$; circles $t=\mbox{\rm const.}$,
$x=\mbox{\rm const.}$ to the second.

\begin{lem}
\label{lem:V}
For every neighbourhood of $\partial V$ in $V$ (or likewise in $B\setminus V$),
there exists a convex torus $T$ inside this neighbourhood, isotopic to
$\partial V$ and satisfying $\# \Gamma_T=2$ and $s(T)=\infty$.
\end{lem}

\begin{proof}
The contact plane $\zeta '$ is spanned by $\frac{\partial}{\partial t}$ and
$\sin (2\pi t)\frac{\partial}{\partial x}+\cos (2\pi t)
\frac{\partial}{\partial y}$. We may choose
$T$ of the form
\[ T=\{ (x,y,t)\in B\co (x,t)\in\gamma\} ,\]
where $\gamma$ is a smooth convex curve in the $(x,t)$--plane, close
to $\partial V\cap\{ y=0\}$. Moreover, we may assume that
$\frac{\partial}{\partial t}$ is
tangent to $\gamma$ only at the two points on $\gamma$ with $t=0$.
The assumption $\varepsilon <1/4$ guarantees that the
singular set of the characteristic foliation $\zeta '|_T$ consists of
the two circles $T\cap\{ t=0\}$. Furthermore, the vector spanning $\zeta '|_T$
away from its singular points always has a non-zero
$\frac{\partial}{\partial x}$--component, and the coefficient functions
of this vector field may be chosen not to depend on the $y$--coordinate.
The two circles $\Gamma =T\cap\{ x=1/2\}$ divide this singular foliation.
Now apply Proposition~\ref{gir2}.
\end{proof}

\begin{proof}[Proof of Proposition~\ref{prop:AA0}] (for $k=1$).
Let $\varphi_0,\varphi$ be as in Lemma~\ref{lem:phi}. Write $\zeta (\varphi_0)$
resp.\ $\zeta (\varphi )$ for the contact structures on $T_{A_0}$ resp.
$T_A$ defined by these functions. Fix a positive real number
$0<\varepsilon < 1/2$. Let $f\co [-1,0]\rightarrow [-1,0]$ be the smooth
function satisfying $\varphi (f(t)) =\varphi_0(t)$ for all $t\in [-1,0]$.
Observe that $f(-1)=-1$ and $f(0)=0$, and $f$ is strictly monotone
increasing.

With $B={\mathbf R}\times ({\mathbf R}/{\mathbf Z})\times {\mathbf R}$ as above,
set
\[ B_{\sigma ,\tau}=\{ (x,y,t)\in B\co \sigma\leq t\leq\tau\} .\]
We continue to write $\zeta (\varphi_0)$, $\zeta (\varphi )$ for the lift
of those contact structures from $T_{A_0}$ resp.\ $T_A$ to~$B$.

Define contact embeddings $F_1$, $F_2$ as follows:
\[ \begin{array}{rrcl}
F_1: & (B_{-1+\varepsilon ,-\varepsilon},\zeta (\varphi_0)) &
\longrightarrow & (B,\zeta (\varphi ))\\
 & (x,y,t) & \longmapsto & (x,y,f(t)),\\
F_2: & (B_{\varepsilon ,1-\varepsilon},\zeta (\varphi_0)) &
\longrightarrow & (B,\zeta (\varphi ))\\
 & (x,y,t) & \longmapsto & (x,x+y,f(t-1)+1). \end{array} \]
Notice that $F_2$ is the composition of contactomorphisms
\[ ({\mathbf x},t)\mapsto (A_0^{-1}{\mathbf x},t-1)\stackrel{F_1}{\longmapsto}
(A_0^{-1}{\mathbf x},f(t-1))\mapsto (AA_0^{-1}{\mathbf x},f(t-1)+1).\]

Fix a positive real number $0<\delta <1/4$. Choose $\varepsilon >0$
sufficiently small such that
\[ -\pi /2 < \varphi_0(-\varepsilon )<\varphi_0(\varepsilon )<\pi /2,\]
\[ -\pi /2 < \varphi (f(-\varepsilon ))<\varphi (f(\varepsilon -1)+1)<\pi /2,\]
and
\[ -\delta <\tan\varphi (f(-\varepsilon ))<\tan\varphi (f(\varepsilon -1)+1)
<\delta.\]
Let $g_1\co [-\varepsilon ,\varepsilon ]\rightarrow {\mathbf R}$ be a
smooth, strictly monotone increasing function such that
\[ g_1(t)=\left\{ \begin{array}{ll} f(t) & \mbox{\rm for}\; -\varepsilon \leq
                     t\leq -\varepsilon /2,\\
        f(t-1)+1 & \mbox{\rm for}\;\; \varepsilon /2 \leq
                     t\leq \varepsilon ,\\
        0 & \mbox{\rm for}\;\; t=0.
\end{array}\right. \]
Let $g_2\co [-\varepsilon ,\varepsilon ]\rightarrow {\mathbf R}$ be a smooth,
monotone increasing function such that
\[ g_2(t)=\left\{ \begin{array}{ll}
0 & \mbox{\rm for} \; -\varepsilon \leq
                     t\leq -\varepsilon /2,\\
1 & \mbox{\rm for} \;\; \varepsilon /2\leq
                     t\leq \varepsilon .\end{array} \right. \]
It is easy to see that $g_2$ can be chosen in such a way that
\[ \int_{-\varepsilon}^{\varepsilon} g_2'(t)\tan\varphi (g_1(t))\, dt =0.\]
For $\eta\in {\mathbf R}$ set
\[ B^{\eta}=B_{-\varepsilon ,\varepsilon}^{\eta}=\{ (x,y,t)\in B\co
x=\eta ,\; -\varepsilon\leq t\leq\varepsilon\} .\]
Define $h_{\eta}\co [-\varepsilon ,\varepsilon ]\rightarrow {\mathbf R}$ by
\[ h_{\eta}(t)=\eta +\eta\int_{-\varepsilon}^t g_2'(t)\tan\varphi (g_1(t))\,
dt .\]
Define
\[ \begin{array}{rrcl}
\psi_{\eta}: & B^{\eta } & \longrightarrow & B\\
 & (\eta ,y,t) & \longmapsto & (h_{\eta }(t),y+g_2(t)\eta ,g_1(t)).
\end{array}\]
Notice that $\psi_{\eta}$ coincides with $F_1$ for $t=-\varepsilon$ and with
$F_2$ for $t=\varepsilon$. Moreover, one easily verifies that $\psi_{\eta}$
is an injective immersion.

We compute
\begin{eqnarray*}
\lefteqn{\psi_{\eta}^*(\cos\varphi (t)\, dx-\sin\varphi (t)\, dy)  =}\\
& = &\cos\varphi (g_1(t))h_{\eta}'\, dt-\sin\varphi (g_1(t))(dy+\eta g_2'(t)
      \, dt)\\
  & = & \sin\varphi (g_1(t))\, dy.
\end{eqnarray*}
It follows that the singular foliation $\psi_{\eta}^{-1}(\zeta
(\varphi )|_{\psi_{\eta}(B^{\eta})})$ is represented by the vector
field $\sin \varphi (g_1(t)) \frac{\partial}{\partial t}$.

The singular foliation $\zeta (\varphi_0)|_{B^{\eta}}$, on the other hand,
is represented by $\sin\varphi_0(t)\frac{\partial}{\partial t}$. We claim
that these two singular foliations are identical as smooth foliations.
Indeed, the two functions $s_1=\sin\varphi (g_1(t))$ and $s_0=\sin\varphi_0(t)$
vanish only at $0\in [-\varepsilon ,\varepsilon ]$ and have positive
derivative there. It follows that either of them can be written as
$s_i=t\cdot \overline{s}_i$ with $\overline{s}_i$ a smooth, nowhere
zero function on $[-\varepsilon ,\varepsilon ]$, so $s_1/s_0$ is smooth and
non-zero on all of $[-\varepsilon ,\varepsilon ]$.

By Proposition~\ref{gir3} there exists a neighbourhood $U$ of $B^{\delta}
\cup B^{1-\delta}$ in $B$ and a contact embedding
\[ F\co (U,\zeta (\varphi_0))\longrightarrow (B,\zeta (\varphi ))\]
that coincides with $F_1$ resp.\ $F_2$ on the common domain of definition,
and with $\psi_{\eta}$ on $B^{\eta}$ for $\eta =\delta$ or $1-\delta$.

By Proposition~\ref{hon} and Lemma~\ref{lem:V}, and with $V$ as in that lemma
(which holds true for the contact structure $\zeta (\varphi_0)$ in place
of~$\zeta '$),
this $F$ extends to a contact embedding
\[ F\co (V_0,\zeta (\varphi_0)):=(B_{-1+\varepsilon ,-\varepsilon}
\cup B_{\varepsilon ,1-\varepsilon}\cup V,\zeta (\varphi_0))\longrightarrow
(B,\zeta (\varphi )).\]

Let $K$ be the Legendrian circle in $T_{A_0}$ defined by
\[ \{ (x,y,t)\in T_{A_0}\co x=t=0\} .\]
Then $F$ induces a contact embedding
\[ (T_{A_0}-\nu K,\zeta (\varphi_0))\longrightarrow (T_A,\zeta (\varphi )),\]
where we may think of the tubular neighbourhood $\nu K$ of $K$ as
\[ \{ (x,y,t)\in T_{A_0}\co -\delta\leq x\leq \delta,\; -\varepsilon\leq t
\leq \varepsilon\} .\]
Again by Proposition~\ref{hon} and Lemma~\ref{lem:V} (adapted
suitably), $(T_A,\zeta (\varphi ))$
is obtained from the manifold $(T_{A_0},\zeta (\varphi_0))$ by
contact $(-1)$--surgery
on~$K$. To verify the sign of this surgery we need to make the following
observations.

Let $\mu$ be a meridian of $\partial\nu K$ defined by $y=0$, say, and let
$\lambda$ be a longitude of $\partial\nu K$ defined by $x=\delta$,
$t=0$. We take $\lambda$ to be oriented in positive $y$--direction,
and $\mu$ to be oriented in counterclockwise direction with respect to
the oriented basis $(\frac{\partial}{\partial t},\frac{\partial}{\partial x})$
of the $(t,x)$--plane. This is consistent with our orientation assumptions
in the definition of contact surgery. Moreover, it is this choice
of longitude that gives $s(\partial\nu K)=\infty$, so the surgery
coefficient $r$ is determined by expressing the attaching map in terms
of $\lambda$ and~$\mu$.

The effect of the map $F$ (up to isotopy) is to send $\lambda$ to
$\lambda$, and $\mu$ to $\mu +\lambda$, as can be checked from
our explicit formulae. So $\mu -\lambda$ maps to $\mu$, which shows that
it is this curve $\mu -\lambda$ on $\partial\nu K$ which becomes
homologically trivial when we glue in a solid torus in place of $\nu K$
to obtain~$T_A$.
\end{proof}

\begin{ack}
This work was done while F.D.\ was visiting the Mathematical
Institute of Leiden University. Part of this stay was supported by
a post-doctoral fellowship from NWO (Netherlands Organisation for
Scientific Research). F.D.\ thanks Leiden University and NWO for their
support.

H.G.\  gratefully acknowledges support by the American Institute of
Mathematics during its programme on Low-Dimensional Contact Geometry.

We also thank Ko Honda for helpful e-mail correspondence.
\end{ack}

\Addresses\recd
\end{document}

%% file: agtout.tex
\def\ifplaintex{\expandafter\ifx\csname documentclass\endcsname\relax}

\def\gtp{{\mathsurround=0pt\it $\cal G\mskip-2mu$eometry \&\ 
$\cal T\!\!$opology $\cal P\!$ublications}}  

\def\recd{{\small Received:\qua\receiveddate\ifx\reviseddate\relax
\else\qquad Revised:\qua\reviseddate\fi\par}} 

\def\volumenumber#1{\def\thevolumenumber{#1}}
\def\volumeyear#1{\def\thevolumeyear{#1}}
\def\papernumber#1{\def\thepapernumber{#1}}
\def\pagenumbers#1#2{\def\startpage{#1}\def\finishpage{#2}}
\def\published#1{\def\publishdate{#1}}

\def\received#1{\def\receiveddate{#1}}
\def\revised#1{\def\reviseddate{#1}}
\def\accepted#1{\def\accepteddate{#1}}

\def\asciiauthors#1{\def\theasciiauthors{#1}}
\def\asciiaddress#1{\def\theasciiaddress{#1}}
\def\asciiemail#1{\def\theasciiemail{#1}}
\def\coverauthors#1{\def\thecoverauthors{#1}}
\long\def\asciiabstract#1{\long\def\theasciiabstract{#1}}

\let\\\par\let\thevolumenumber\relax
\let\thepapernumber\relax\let\thevolumeyear\relax\let\startpage\relax
\let\finishpage\relax\let\publishdate\relax\let\receiveddate\relax
\let\reviseddate\relax\let\accepteddate\relax\let\theasciititle\relax
\let\theasciiauthors\relax\let\theasciiaddress\relax
\let\theasciiabstract\relax

\let\thecoverauthors\relax\let\theasciiemail\relax

\ifplaintex
\font\logobig=cmssbx10 scaled 3836
\font\logomed=cmssbx10 scaled 2557
\else
\font\logobig=cmssbx10 scaled 4200
\font\logomed=cmssbx10 scaled 2800
\fi

\long\def\makeagttitle{   
\count0=\startpage
\agt\hfill      
\hbox to 45truept{\vbox to 0pt{\vglue -13truept{\logomed A\kern -.37em{\logobig 
T}\kern -.38em G}\vss}\hss}
\break
{\small Volume \thevolumenumber\ (\thevolumeyear)
\startpage--\finishpage\nl
Published: \publishdate}

\vglue .25truein

{\parskip=0pt\leftskip 0pt plus
1fil\def\\{\par\smallskip}{\Large\bf\thetitle}\par\medskip} \vglue
0.05truein

{\parskip=0pt\leftskip 0pt plus 1fil\def\\{\par}{\sc\theauthors}
\par\medskip}%
 
\vglue 0.03truein 

{\small\leftskip 25truept\rightskip 25truept{\bf Abstract}\stdspace\theabstract

{\bf AMS Classification}\stdspace\theprimaryclass
\ifx\thesecondaryclass\relax\else; \thesecondaryclass\fi\par
{\bf Keywords}\stdspace \thekeywords\par}\vglue 7truept

}   

\ifplaintex
\font\phead=cmsl9 scaled 950
\font\pnum=cmbx10 scaled 913
\font\pfoot=cmsl9 scaled 950
\headline{\vbox to 0pt{\vskip -4.5mm\line{\small\phead\ifnum
\count0=\startpage ISSN 1472-2739 (on-line) 1472-2747 (printed)
\hfill {\pnum\folio}\else\ifodd\count0\def\\{ }%
\ifx\theshorttitle\relax\thetitle\else\theshorttitle\fi\hfill{\pnum\folio}
\else\def\\{ and }{\pnum\folio}\hfill\ifx\theshortauthors\relax\theauthors
\else\theshortauthors\fi\fi\fi}\vss}}
\footline{\vbox to 0pt{\vglue 0mm\line{\small\pfoot\ifnum\count0=\startpage
\copyright\ \gtp\hfill\else
\agt, Volume \thevolumenumber\ (\thevolumeyear)\hfill\fi}\vss}}
\else
\font=cmsl9 scaled 1050
\font\lnum=cmbx10 
\font=cmsl9 scaled 1050
\makeatletter
\def\@oddhead{{\small\lhead\ifnum\count0=\startpage ISSN 1472-2739 
(on-line) 1472-2747 (printed)\hfill {\lnum\number\count0}\else\ifodd\count0
\def\\{ }\ifx\theshorttitle\relax \thetitle \else\theshorttitle\fi\hfill
{\lnum\number\count0}\else\def\\{ and }{\lnum\number\count0}
\hfill\ifx\theshortauthors\relax 
\theauthors\else\theshortauthors\fi\fi\fi}}\def\@evenhead{\@oddhead}
\def\@oddfoot{\small\lfoot\ifnum\count0=\startpage\copyright\ \gtp\hfill\else
\agt, Volume \thevolumenumber\ (\thevolumeyear)\hfill\fi}
\def\@evenfoot{\@oddfoot}
\makeatother
\fi
\let\maketitlepage\makeagttitle
\let\makeshorttitle\maketitlepage
\let\maketitle\maketitlepage

\newwrite\gtoutfile
\long\gdef\makeheadfile{  
{\def\\{, }\def\s{ }
\immediate\openout\gtoutfile head.xxx
\immediate\write\gtoutfile{To: math@arxiv.org}
\immediate\write\gtoutfile{Subject: put OR rep NNNNN:ppppp}
\immediate\write\gtoutfile{--text follows this line--}
\immediate\write\gtoutfile{Proxy-for: \ifx\theasciiauthors\relax
\theauthors\else\theasciiauthors\fi\s<\ifx\theasciiemail\relax\theemail\else\theasciiemail\fi>}
\immediate\write\gtoutfile{\noexpand\\}
\immediate\write\gtoutfile{Authors: \ifx\theasciiauthors\relax
\theauthors\else\theasciiauthors\fi}
{\def\\{ }\immediate\write\gtoutfile{Title: \ifx\theasciititle\relax
\thetitle\else\theasciititle\fi}}
\immediate\write\gtoutfile{Subj-class: GT or SG, GR etc}
\immediate\write\gtoutfile{MSC-class: \theprimaryclass\ifx\thesecondaryclass\relax\else, \thesecondaryclass\fi}
\immediate\write\gtoutfile{Journal-ref: Algebr. Geom. Topol. \thevolumenumber\s
(\thevolumeyear) \startpage-\finishpage}
\immediate\write\gtoutfile{Comments: Published by Algebraic and
Geometric Topology at}
\immediate\write\gtoutfile{\s\s\s  http://www.maths.warwick.ac.uk/agt/AGTVol\thevolumenumber/agt-\thevolumenumber-\thepapernumber.abs.html}
\immediate\write\gtoutfile{\noexpand\\}
\immediate\write\gtoutfile{}
\ifx\theasciiabstract\relax
\immediate\write\gtoutfile{\theabstract}\else
\immediate\write\gtoutfile{\theasciiabstract}\fi
\immediate\write\gtoutfile{}
\immediate\write\gtoutfile{\noexpand\\}
\immediate\write\gtoutfile{}
\immediate\closeout\gtoutfile}}  

\def\maketitlepage{\makeagttitle\makeheadfile}
\let\makeshorttitle\maketitlepage
\let\maketitle\maketitlepage

\def\ifplaintex{\expandafter\ifx\csname documentclass\endcsname\relax}

\def\gtp{{\mathsurround=0pt\it $\cal G\mskip-2mu$eometry \&\ 
$\cal T\!\!$opology $\cal P\!$ublications}}  

\def\recd{{\small Received:\qua\receiveddate\ifx\reviseddate\relax
\else\qquad Revised:\qua\reviseddate\fi\par}} 

\def\volumenumber#1{\def\thevolumenumber{#1}}
\def\volumeyear#1{\def\thevolumeyear{#1}}
\def\papernumber#1{\def\thepapernumber{#1}}
\def\pagenumbers#1#2{\def\startpage{#1}\def\finishpage{#2}}
\def\published#1{\def\publishdate{#1}}

\def\received#1{\def\receiveddate{#1}}
\def\revised#1{\def\reviseddate{#1}}
\def\accepted#1{\def\accepteddate{#1}}

\def\asciiauthors#1{\def\theasciiauthors{#1}}
\def\asciiaddress#1{\def\theasciiaddress{#1}}
\def\asciiemail#1{\def\theasciiemail{#1}}
\def\coverauthors#1{\def\thecoverauthors{#1}}
\long\def\asciiabstract#1{\long\def\theasciiabstract{#1}}

\let\\\par\let\thevolumenumber\relax
\let\thepapernumber\relax\let\thevolumeyear\relax\let\startpage\relax
\let\finishpage\relax\let\publishdate\relax\let\receiveddate\relax
\let\reviseddate\relax\let\accepteddate\relax\let\theasciititle\relax
\let\theasciiauthors\relax\let\theasciiaddress\relax
\let\theasciiabstract\relax

\let\thecoverauthors\relax\let\theasciiemail\relax

\ifplaintex
\font\logobig=cmssbx10 scaled 3836
\font\logomed=cmssbx10 scaled 2557
\else
\font\logobig=cmssbx10 scaled 4200
\font\logomed=cmssbx10 scaled 2800
\fi

\long\def\makeagttitle{   
\count0=\startpage
\agt\hfill      
\hbox to 45truept{\vbox to 0pt{\vglue -13truept{\logomed A\kern -.37em{\logobig 
T}\kern -.38em G}\vss}\hss}
\break
{\small Volume \thevolumenumber\ (\thevolumeyear)
\startpage--\finishpage\nl
Published: \publishdate}

\vglue .25truein

{\parskip=0pt\leftskip 0pt plus
1fil\def\\{\par\smallskip}{\Large\bf\thetitle}\par\medskip} \vglue
0.05truein

{\parskip=0pt\leftskip 0pt plus 1fil\def\\{\par}{\sc\theauthors}
\par\medskip}%
 
\vglue 0.03truein 

{\small\leftskip 25truept\rightskip 25truept{\bf Abstract}\stdspace\theabstract

{\bf AMS Classification}\stdspace\theprimaryclass
\ifx\thesecondaryclass\relax\else; \thesecondaryclass\fi\par
{\bf Keywords}\stdspace \thekeywords\par}\vglue 7truept

}   

\ifplaintex
\font\phead=cmsl9 scaled 950
\font\pnum=cmbx10 scaled 913
\font\pfoot=cmsl9 scaled 950
\headline{\vbox to 0pt{\vskip -4.5mm\line{\small\phead\ifnum
\count0=\startpage ISSN 1472-2739 (on-line) 1472-2747 (printed)
\hfill {\pnum\folio}\else\ifodd\count0\def\\{ }%
\ifx\theshorttitle\relax\thetitle\else\theshorttitle\fi\hfill{\pnum\folio}
\else\def\\{ and }{\pnum\folio}\hfill\ifx\theshortauthors\relax\theauthors
\else\theshortauthors\fi\fi\fi}\vss}}
\footline{\vbox to 0pt{\vglue 0mm\line{\small\pfoot\ifnum\count0=\startpage
\copyright\ \gtp\hfill\else
\agt, Volume \thevolumenumber\ (\thevolumeyear)\hfill\fi}\vss}}
\else
\font=cmsl9 scaled 1050
\font\lnum=cmbx10 
\font=cmsl9 scaled 1050
\makeatletter
\def\@oddhead{{\small\lhead\ifnum\count0=\startpage ISSN 1472-2739 
(on-line) 1472-2747 (printed)\hfill {\lnum\number\count0}\else\ifodd\count0
\def\\{ }\ifx\theshorttitle\relax \thetitle \else\theshorttitle\fi\hfill
{\lnum\number\count0}\else\def\\{ and }{\lnum\number\count0}
\hfill\ifx\theshortauthors\relax 
\theauthors\else\theshortauthors\fi\fi\fi}}\def\@evenhead{\@oddhead}
\def\@oddfoot{\small\lfoot\ifnum\count0=\startpage\copyright\ \gtp\hfill\else
\agt, Volume \thevolumenumber\ (\thevolumeyear)\hfill\fi}
\def\@evenfoot{\@oddfoot}
\makeatother
\fi
\let\maketitlepage\makeagttitle
\let\makeshorttitle\maketitlepage
\let\maketitle\maketitlepage

\newwrite\gtoutfile
\long\gdef\makeheadfile{  
{\def\\{, }\def\s{ }
\immediate\openout\gtoutfile head.xxx
\immediate\write\gtoutfile{To: math@arxiv.org}
\immediate\write\gtoutfile{Subject: put OR rep NNNNN:ppppp}
\immediate\write\gtoutfile{--text follows this line--}
\immediate\write\gtoutfile{Proxy-for: \ifx\theasciiauthors\relax
\theauthors\else\theasciiauthors\fi\s<\ifx\theasciiemail\relax\theemail\else\theasciiemail\fi>}
\immediate\write\gtoutfile{\noexpand\\}
\immediate\write\gtoutfile{Authors: \ifx\theasciiauthors\relax
\theauthors\else\theasciiauthors\fi}
{\def\\{ }\immediate\write\gtoutfile{Title: \ifx\theasciititle\relax
\thetitle\else\theasciititle\fi}}
\immediate\write\gtoutfile{Subj-class: GT or SG, GR etc}
\immediate\write\gtoutfile{MSC-class: \theprimaryclass\ifx\thesecondaryclass\relax\else, \thesecondaryclass\fi}
\immediate\write\gtoutfile{Journal-ref: Algebr. Geom. Topol. \thevolumenumber\s
(\thevolumeyear) \startpage-\finishpage}
\immediate\write\gtoutfile{Comments: Published by Algebraic and
Geometric Topology at}
\immediate\write\gtoutfile{\s\s\s  http://www.maths.warwick.ac.uk/agt/AGTVol\thevolumenumber/agt-\thevolumenumber-\thepapernumber.abs.html}
\immediate\write\gtoutfile{\noexpand\\}
\immediate\write\gtoutfile{}
\ifx\theasciiabstract\relax
\immediate\write\gtoutfile{\theabstract}\else
\immediate\write\gtoutfile{\theasciiabstract}\fi
\immediate\write\gtoutfile{}
\immediate\write\gtoutfile{\noexpand\\}
\immediate\write\gtoutfile{}
\immediate\closeout\gtoutfile}}  

\def\maketitlepage{\makeagttitle\makeheadfile}
\let\makeshorttitle\maketitlepage
\let\maketitle\maketitlepage